\documentclass[12pt]{article}
\usepackage{amsfonts}

\usepackage{amsmath}
\usepackage[]{graphicx}

\newtheorem{theorem}{Theorem}

\newtheorem{definition}[theorem]{Definition}
\newtheorem{example}[theorem]{Example}
\newtheorem{lemma}[theorem]{Lemma}

\newtheorem{remark}[theorem]{Remark}
\newenvironment{proof}[1][Proof]{\textbf{#1.} }{\
\rule{0.5em}{0.5em}}

\title{LULU operators for functions of continuous argument}   
\author{Roumen Anguelov\\
Department of Mathematics and Applied Mathematics\\
University of Pretoria, South
Africa\\
roumen.anguelov@up.ac.za\vspace{5mm}\\
Carl Rohwer\\
Department of Mathematical Sciences\\
University of Stellenbosch, South Africa\\chr@sun.ac.za}
\date{}

\begin{document}
\maketitle
\begin{abstract}
The LULU operators, well known in the nonlinear multiresolution
analysis of sequences, are extended to functions defined on a
continuous domain, namely, a real interval. We show that the
extended operators replicate the essential properties of their
discrete counterparts. More precisely, they form a fully ordered
semi-group of four elements, preserve the local trend and the
total variation.
\end{abstract}

\section{Introduction}

The well known (linear-) Functional Analysis fits in appropriately
in the theory of linear smoothers. Typically a smoother is
designed to pass sequences that are samplings of functions with
low frequencies with minimal distortion (error), but to map high
frequencies near to the zero sequence. The least squares norm is
appropriate, for various other reasons, but also that such
smoothers then map a element $x_i$ onto weighted averages of
sequence elements in a ``window'' around $x_i$, say $\{x_{i-n},
\ldots, x_i, \ldots, x_{i+n} \}$. The linearity also ensures that
sequences of elements that are generated identically,
independently distributed from a very general symmetrical
distribution $e$, are rapidly mapped near zero due to the Central
Limit Theorem. When we want to smooth a sequence, we can choose to
construct a convenient smoother that is a "bandpass" filter,and
practically remove high frequencies. The design of such filters is
done in a well established theory of Digital Filters. The book by
Hamming \cite{Hamming} is well known and instructive. The
essential background theory is Fourier Analysis. This works well
because the basic trigonometric functions $\sin$ and $\cos$ are
eigensequences of linear operators, and the ``transfer function''
approximates with eigenvalues near one in the frequencies that are
to pass, and eigenvalues near $0$ where the frequencies are to be
reduced to near zero.

Low pass filters should therefore marginally distort sequences
that are samplings of functions that have Fourier expansions that
converge fast. This is well known to be related to continuity of
lower derivatives.   Discontinuities in low order derivatives
result in slowly converging Fourier expansions and the digital
filters will remove the high frequencies, distorting significantly
the sequence of samplings. A typical bad case is isolated
impulsive noise added spuriously. This necessitates the
presmoothing by nonlinear Smoothers, of which the median
smoothers, popularised by Tukey, are well known. Since
eigensequence analysis is not natural, nor easily justifiable, for
nonlinear operators the lack of a theory for analysis and design
was considered to be difficult, if not impossible \cite{Serra}.
Design was generally essentially considered to be an art.

Over the last twenty years a theory for Nonlinear (general)
Smoothers, that is based on order structure and min/max
operations, has been developed and demonstrated to be very
consistent and useful, even able to explain most of the ``good''
behaviour of the (related) median smoothers, as well as their
``enigmatic'' behaviour. A monogram presenting the so-called
$LULU$-theory and its motivation and development appeared in 2005
\cite{Rohwerbook}. The theory is based on compositions of two
types of smoothers $L_n$ and $U_n$. They are Morphological Filters
with special properties.

One of the powerful ideas resulting from this theory, was the
development of Nonlinear Multiresolution Analysis. This was done
using the heuristic ideas from Fourier Analysis and Wavelet
Analysis. It resulted eventually in Discrete Pulse Transforms
\cite{Laurie}. These transforms may turn out to be as useful for
vision as the Fourier Transforms are for hearing
\cite{RohwerWild}.

When applications of Wavelet Transforms (and Fourier Transforms)
are under discussion it is natural for understanding to consider
samplings of ``band limited'' functions as ideal candidates for
such decompositions, both for theoretical derivation and practical
applications. For Nonlinear Decomposition there has been a lack of
such a relation between the theory of real functions and the
theory of the sequences that are samplings of these functions.
Generalising $LULU$-operators and the associated theory and
concepts to functions is the first appropriate attempt towards
establishing such a link.

Central in the $LULU$-theory for sequences are the class of
locally monotone sequences ${\cal M}_n$ defined as the sets of all
sequences $x$ that have $\{x_i, x_{i+1}, \ldots, x_{i+n+1} \}$
monotone for each index $i$. We need to establish natural links
between these classes and classes of real functions of which the
sequences can be considered as samplings. Also natural to
$LULU$-theory is the Total Variation as norm \cite{Rohwerbook}.
There is a clue to establishing links with standard Real Analysis,
as is typically presented by Royden in the first few chapters of
his book \cite{Royden}. Total Variation of functions and local
monotonicity are linked to the derivative in this theory. We seek
to extend and solidify these links with established Real Analysis
of functions. To do this we, look at the basic ideas of the
$LULU$-theory for sequences \cite{Rohwerbook}. We start directly
with the definitions of the ``atoms'' involved.

Given a bi-infinite sequence $\xi=(\xi_i)_{i\in\mathbb{Z}}$ and
$n\in\mathbb{N}$ the operators $L_n$ and $U_n$ are defined as
follows
\begin{eqnarray*}
(L_n
\xi)_i=\max\{\min\{\xi_{i-n},...,\xi_i\},...,\min\{\xi_{i},...,\xi_{i+n}\}\},\
i\in\mathbb{N}\\
(U_n
\xi)_i=\min\{\max\{\xi_{i-n},...,\xi_i\},...,\max\{\xi_{i},...,\xi_{i+n}\}\},\
i\in\mathbb{N}
\end{eqnarray*}
In analogy with the above discrete LULU operators, for a given
$\delta>0$ the basic smoothers $L_\delta$ and $U_\delta$ in the
LULU theory are defined for functions on $\Omega$ through the
concepts of the so called lower and upper $\delta$-envelopes of
these functions. These definitions are given in Section 2, where
it is also shown that the operators $L_\delta$ and $U_\delta$
preserve essential properties of their discrete counterparts. In
particular, the operators $L_\delta$ and $U_\delta$ generate
through composition a fully ordered four element semi-group, also
called a strong LULU structure, as opposed to the general weak
$LULU$-structure of recursions and dilations, with a 6 element
semi-group that is only partially ordered. In particular, the
crucial inequality $L_nU_n \geq U_nL_n$ holds in a strong
structure. This issue is dealt with in Section 3. In Sections 4
and 5 we discuss the preservation of the trend and the total
variation respectively.

\newpage
\section{The basic smoothers $L_\delta$ and $U_\delta$}

Let $\mathcal{A}(\Omega)$ denote the set of all bounded real
functions defined on a real interval $\Omega\subseteq\mathbb{R}$.
Let $B_{\delta}(x)$ denote the closed $\delta$-neighborhood of $x$
in $\Omega $, that is, $B_{\delta }(x)=\{y\in \Omega
:|x-y|\leq\delta \}$.
The pair of mappings $I$, $S:\mathcal{A}(\Omega )\rightarrow \mathcal{A}%
(\Omega )$ defined by
\begin{eqnarray}
I(f)(x) &=&\sup_{\delta >0}\inf \{f(y):y\in B_{\delta }(x)\}, \
x\in\Omega,
\label{lbf} \\
S(f)(x) &=&\inf_{\delta >0}\sup \{f(y):y\in B_{\delta }(x)\}, \
x\in\Omega, \label{ubf}
\end{eqnarray}%
are called lower Baire, and upper Baire operators, respectively,
\cite{Sendov}. We consider on $\mathcal{A}(\Omega)$ the point-wise
defined partial order, that is, for any
$f,g\in\mathcal{A}(\Omega)$
\begin{equation}\label{forder}
f\leq g\Longleftrightarrow f(x)\leq g(x),\ x\in\Omega.
\end{equation}
Then the lower and upper Baire operators can be defined in the
following equivalent way. For every $f\in \mathcal{A}(\Omega )$
the function $I(f)$ is the maximal lower semi-continuous function
which is not greater than $f$. Hence, it is also called lower
semi-continuous envelope. In a similar way, $S(f)$ is the smallest
upper semi-continuous function which is not less than $f$ and is
called the upper semi-continuous envelope of $f$. In analogy with
$I(f)$ and $S(f)$ we call the functions
\begin{eqnarray}
I_{\delta }(f)(x) &=&\inf \{f(y):y\in B_{\delta }(x)\},\
x\in\Omega,
\label{ldenvelop} \\
S_{\delta }(f)(x) &=&\sup \{f(y):y\in B_{\delta }(x)\},\
x\in\Omega, \label{udenvelop}
\end{eqnarray}%
a lower $\delta$-envelope of $f$ and an upper $\delta $-envelope
of $f$, respectively.

It is easy to see from (\ref{ldenvelop}) and (\ref{udenvelop})
that for every $\delta_1,\delta_2>0$
\begin{equation}
I_{\delta _{1}}\circ I_{\delta _{2}}=I_{\delta _{1}+\delta _{2}}\
\ ,\ \ \  S_{\delta _{1}}\circ S_{\delta _{2}}=S_{\delta
_{1}+\delta _{2}}\label{IScomp}
\end{equation}
Furthermore, the operators $I_\delta$ and $S_\delta$, $\delta>0$,
as well as $I$ and $S$ are all monotone increasing with respect to
the order (\ref{forder}), that is, for every
$f,g\in\mathcal{A}(\Omega)$
\begin{equation}\label{IdSdmon}f\leq g\Longrightarrow I_\delta(f)\leq
I_\delta(g),\ S_\delta(f)\leq S_\delta(g),\ I(f)\leq I(g),\
S(f)\leq S(g).
\end{equation}

The following operators can be considered as continuous analogues
of the discrete $LULU$ operators given in the Introduction:
\begin{equation}\label{defLULU}
L_{\delta }=S_{\frac{\delta}{2} }\circ I_{\frac{\delta}{2} }\ ,\ \
U_{\delta }=I_{\frac{\delta}{2} }\circ S_{\frac{\delta}{2}}\ .
\end{equation}
We will show that these operators have similar properties to their
discrete counterparts. Let us note that they inherit the
monotonicity with respect to the functional argument from the
operators $I_\delta$ and $S_\delta$, see (\ref{IdSdmon}), that is,
for $f,g\in\mathcal{A}(\Omega)$
\begin{equation}\label{LdUdfmon}f\leq g\Longrightarrow L_\delta(f)\leq
L_\delta(g),\ U_\delta(f)\leq U_\delta(g).
\end{equation}

\begin{theorem}\label{tLdfleqf}
For every $f\!\in\!\mathcal{A}(\Omega)$ and $\delta\!>\!0$ we have
$L_\delta(f)\leq f$, $U_\delta(f)\geq f$.
\end{theorem}
\begin{proof}
Let $f\in\mathcal{A}(\Omega)$, $\delta>0$. It follows from the
definition of $I_\delta$ that for any $x\in\Omega$ we have
$I_\frac{\delta}{2}(f)(y)\leq f(x),\ y\in B_\frac{\delta}{2}(x)$.
Therefore
\[
L_\delta(f)(x)=S_\frac{\delta}{2}(I_\frac{\delta}{2}(f))(x)=\sup\{I_\frac{\delta}{2}(f)(y):y\in
B_\frac{\delta}{2}(x)\}\leq f(x),\  x\in\Omega.
\]
The second inequality in the theorem is proved in a similar way.
\end{proof}

\begin{theorem}\label{tLdUdmon}
The operator $L_\delta$ is monotone decreasing on $\delta$ while
the operator $U_\delta$ is monotone increasing on $\delta$, that
is, for any $f\in\mathcal{A}(\Omega)$ and $0<\delta_1\leq\delta_2$
we have $L_{\delta_1}(f)\geq L_{\delta_2}(f)$,
$U_{\delta_1}(f)\leq U_{\delta_2}(f)$.
\end{theorem}
\begin{proof}
Let $\delta_2>\delta_1>0$. Using the properties (\ref{IScomp}) the
operator $L_{\delta_2}$ can be represented in the form
\[%
L_{\delta_2}=S_\frac{\delta_2}{2}\circ
I_\frac{\delta_2}{2}=S_\frac{\delta_1}{2}\circ
S_{\frac{\delta_2-\delta_1}{2}}\circ
I_{\frac{\delta_2-\delta_1}{2}} \circ
I_\frac{\delta_1}{2}=S_\frac{\delta_1}{2}\circ
L_{\delta_2-\delta_1} \circ I_\frac{\delta_1}{2}.%
\]
It follows from Theorem \ref{tLdfleqf} that for every
$f\in\mathcal{A}(\Omega)$ we have
$L_{\delta_2-\delta_1}(I_\frac{\delta_1}{2}(f))\leq
I_\frac{\delta_1}{2}(f)$. Hence using the monotonicity of the
operator $S_\delta$ given in (\ref{IdSdmon}) we obtain
\[%
L_{\delta_2}(f)=S_\frac{\delta_1}{2}(L_{\delta_2-\delta_1}(I_\frac{\delta_1}{2}(f)))\leq
S_\frac{\delta_1}{2}(I_\frac{\delta_1}{2}(f))=L_{\delta_1}(f),\
f\in\mathcal{A}(\Omega).%
\]
The inequality $U_{\delta_1}(f)\geq U_{\delta_2}(f)$ is proved in
a similar way.
\end{proof}

The next lemma is useful in dealing with compositions of
$I_\delta$ and $S_\delta$.
\begin{lemma}\label{tIdSdId} We have
$I_\delta\circ S_\delta\circ I_\delta=I_\delta$, $S_\delta\circ
I_\delta\circ S_\delta=S_\delta$.
\end{lemma}
\begin{proof}
Using the monotonicity of $I_\delta$, see (\ref{IdSdmon}), and
Theorem \ref{tLdfleqf} for $f\in \mathcal{A}(\Omega)$ we have
$(I_\delta\circ S_\delta\circ I_\delta)(f)=I_\delta
(L_{2\delta}(f))\leq I_\delta (f)$. On the other hand, applying
Theorem \ref{tLdfleqf} to $U_{2\delta}$ we obtain $(I_\delta\circ
S_\delta\circ I_\delta)(f)=U_{2\delta}(I_\delta(f))\geq
I_\delta(f)$. Therefore $(I_\delta\circ S_\delta\circ
I_\delta)(f)=I_\delta(f)$, $f\in\mathcal{A}(\Omega)$. The second
equality is proved similarly.
\end{proof}

\begin{theorem}\label{tLdUdabsorb}
For every $\delta_1,\delta_2>0$ we have $L_{\delta_1 }\circ
L_{\delta_2 } =L_{\max\{\delta_1,\delta_2\}}$ and $U_{\delta_1
}\circ U_{\delta_2 }=U_{\max\{\delta_1,\delta_2\}}$.
\end{theorem}
\begin{proof}
We will only prove the first equality since the proof of the
second one is done in a similar manner. Let first
$\delta_2>\delta_1>0$. Using property (\ref{IScomp}) and Lemma
\ref{tIdSdId} we obtain
\begin{eqnarray*}
L_{\delta_1}\circ L_{\delta_2}&=&(S_\frac{\delta_1}{2}\circ
I_\frac{\delta_1}{2})\circ(S_\frac{\delta_2}{2}\circ
I_\frac{\delta_2}{2})\ =\ (S_\frac{\delta_1}{2}\circ
I_\frac{\delta_1}{2}\circ S_\frac{\delta_1}{2})\circ
(S_{\frac{\delta_2-\delta_1}{2}}\circ I_\frac{\delta_2}{2})\\
&=& S_\frac{\delta_1}{2}\circ S_{\frac{\delta_2-\delta_1}{2}}\circ
I_\frac{\delta_2}{2}\ =\ S_\frac{\delta_2}{2}\circ
I_\frac{\delta_2}{2}\ =\ L_{\delta_2}.
\end{eqnarray*}
If $\delta_1>\delta_2>0$ in a similar way we have
\begin{eqnarray*}
L_{\delta_1}\circ L_{\delta_2}&=&(S_{\delta_1}\circ
I_{\delta_1})\circ(S_{\delta_2}\circ I_{\delta_2})\ =\
(S_{\delta_1}\circ I_{\delta_1-\delta_2})\circ (I_{\delta_2}\circ
S_{\delta_2}\circ I_{\delta_2})\\
&=&S_{\delta_1}\circ I_{\delta_1-\delta_2}\circ I_{\delta_2}\ =\
S_{\delta_1}\circ I_{\delta_1}\ =\ L_{\delta_1}.
\end{eqnarray*}
The proof in the case when $\delta_2\!=\!\delta_1\!>\!0$ follows
from either of the above identities where $S_{\delta_2-\delta_1}$
or $I_{\delta_1-\delta_2}$ respectively are  replaced by the
identity operator.
\end{proof}

Important properties of smoothing operators are their idempotence
and co-idempotence. Hence the significance of the next theorem.
\begin{theorem}\label{tLdUdidemp}
The operators $L_\delta$ and $U_\delta$ are both idempotent and
co-idempotent, that is, $L_{\delta }\circ L_{\delta }
=L_{\delta}$, $U_{\delta }\circ U_{\delta } =U_{\delta}$,
$(id-L_{\delta })\circ (id-L_{\delta })=id-L_{\delta }$,
$(id-U_{\delta })\circ (id-U_{\delta })=id-U_{\delta }$, where
$id$ denotes the identity operator.
\end{theorem}
\begin{proof}
The idempotence of $L_\delta$ and $U_\delta$ follows directly from
Theorem \ref{tLdUdabsorb}. The co-idempotence of the operator
$L_\delta$ is equivalent to $L_\delta\circ(id-L_\delta)=0$. Using
the first inequality in Theorem \ref{tLdfleqf} one can easily
obtain $L_\delta\circ(id-L_\delta)\geq 0$. Hence, for the
co-idempotence of $L_\delta$ it remains to show that
$L_\delta\circ(id-L_\delta)\leq 0$. Assume the opposite. Namely,
there exists a function $f\in\mathcal{A}(\Omega)$ and $x\in\Omega$
such that $(L_\delta\circ(id-L_\delta))(f)(x)>0$. Let
$\varepsilon>0$ be such that
$(L_\delta\circ(id-L_\delta))(f)(x)>\varepsilon>0$. Using the
definition of $L_\delta$ the above inequality implies that there
exists $y\in B_\frac{\delta}{2}(x)$ such that for every $z\in
B_\frac{\delta}{2}(y)$ we have $(id-L_\delta)(f)(z)>\varepsilon$,
or equivalently %
\begin{equation}\label{ineq3tidemp}
f(z)>L_\delta(f)(z)+\varepsilon,\ z\in B_\frac{\delta}{2}(y).
\end{equation}
For every $z\in B_\frac{\delta}{2}(y)$ we also have
$L_\delta(f)(z)\geq I_\frac{\delta}{2}(f)(y)=\inf\{f(t):t\in
B_\frac{\delta}{2}(y)\}$. Hence there exists $t\in
B_\frac{\delta}{2}(y)$ such that
$f(t)<I_\frac{\delta}{2}(f)(y)+\varepsilon\leq
L_\delta(f)(z)+\varepsilon$, $z\in B_\frac{\delta}{2}(y)$. Taking
$z=t$ in the above inequality we obtain
$f(t)<L_\delta(f)(t)+\varepsilon$, which contradicts
(\ref{ineq3tidemp}). The co-idempotence of $U_\delta$ is proved in
a similar way.
\end{proof}

\begin{example} The figures below illustrate graphically the smoothing
effect of the operators $L_\delta$, $U_\delta$ and their
compositions. The graph of function $f$ is given by dotted lines.

\begin{center}
\includegraphics[scale=0.3]{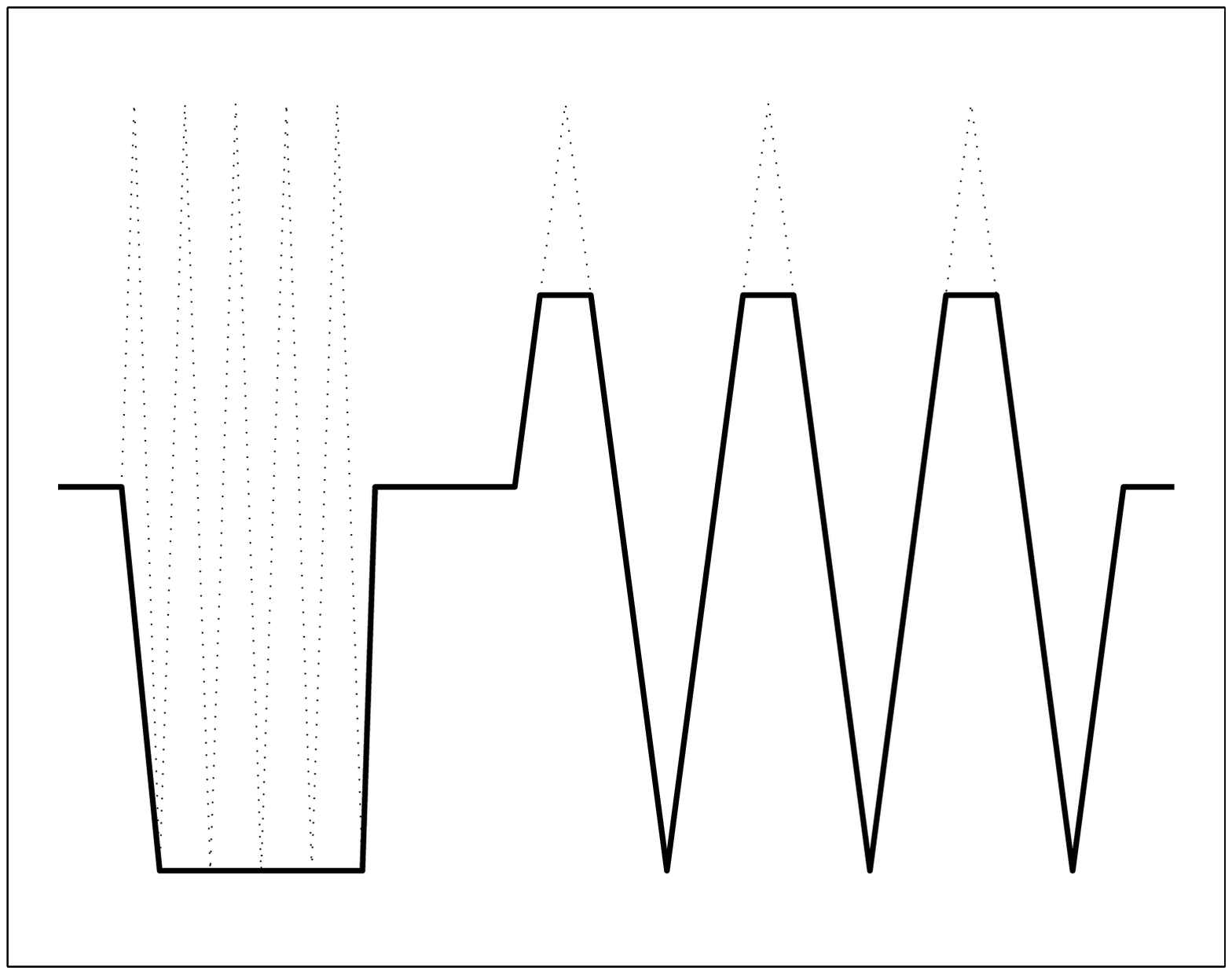}\hspace{5mm}
\includegraphics[scale=0.3]{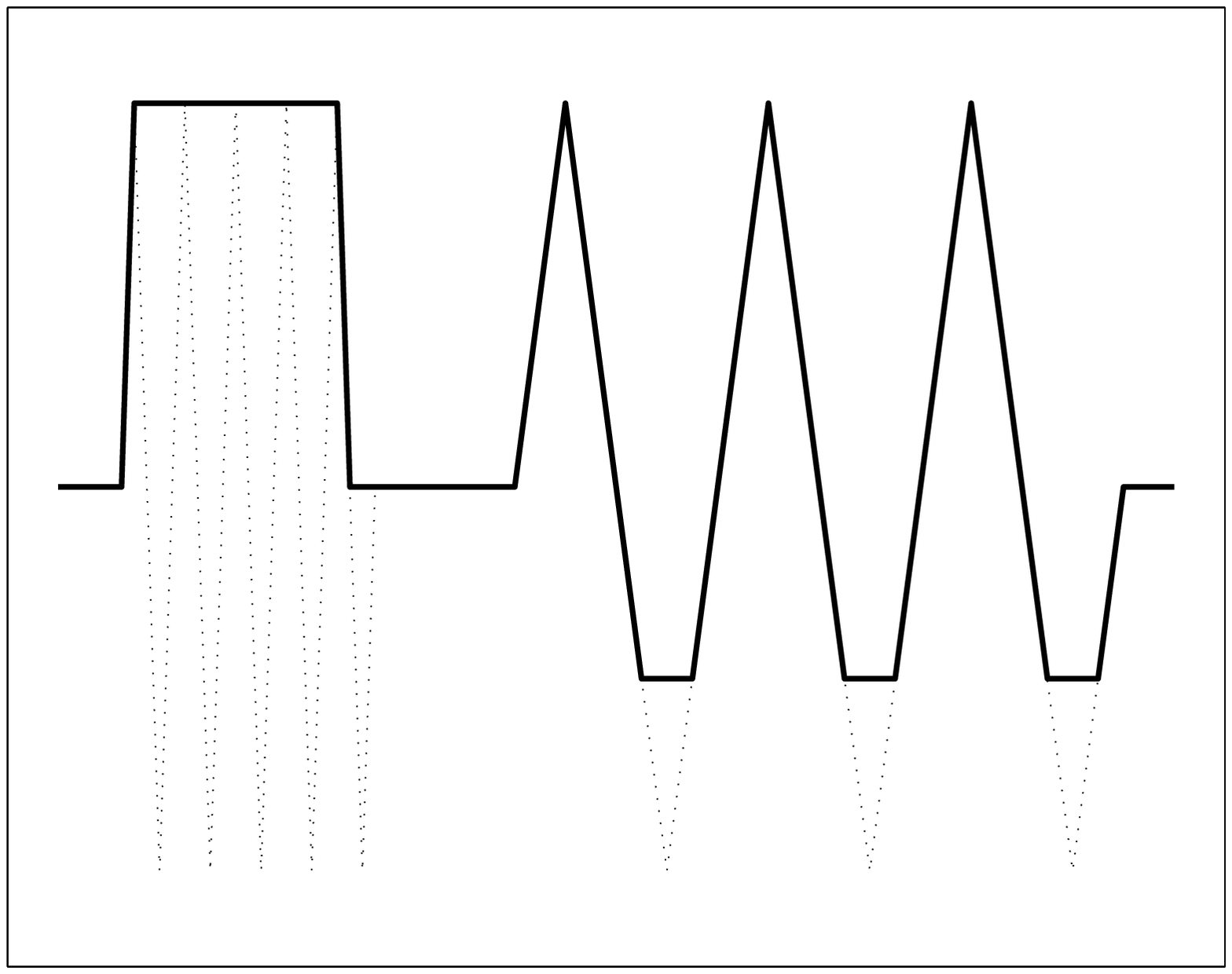}\\\nopagebreak
The functions $L_\delta(f)$ and $U_\delta(f)$\\[6pt]
\includegraphics[scale=0.3]{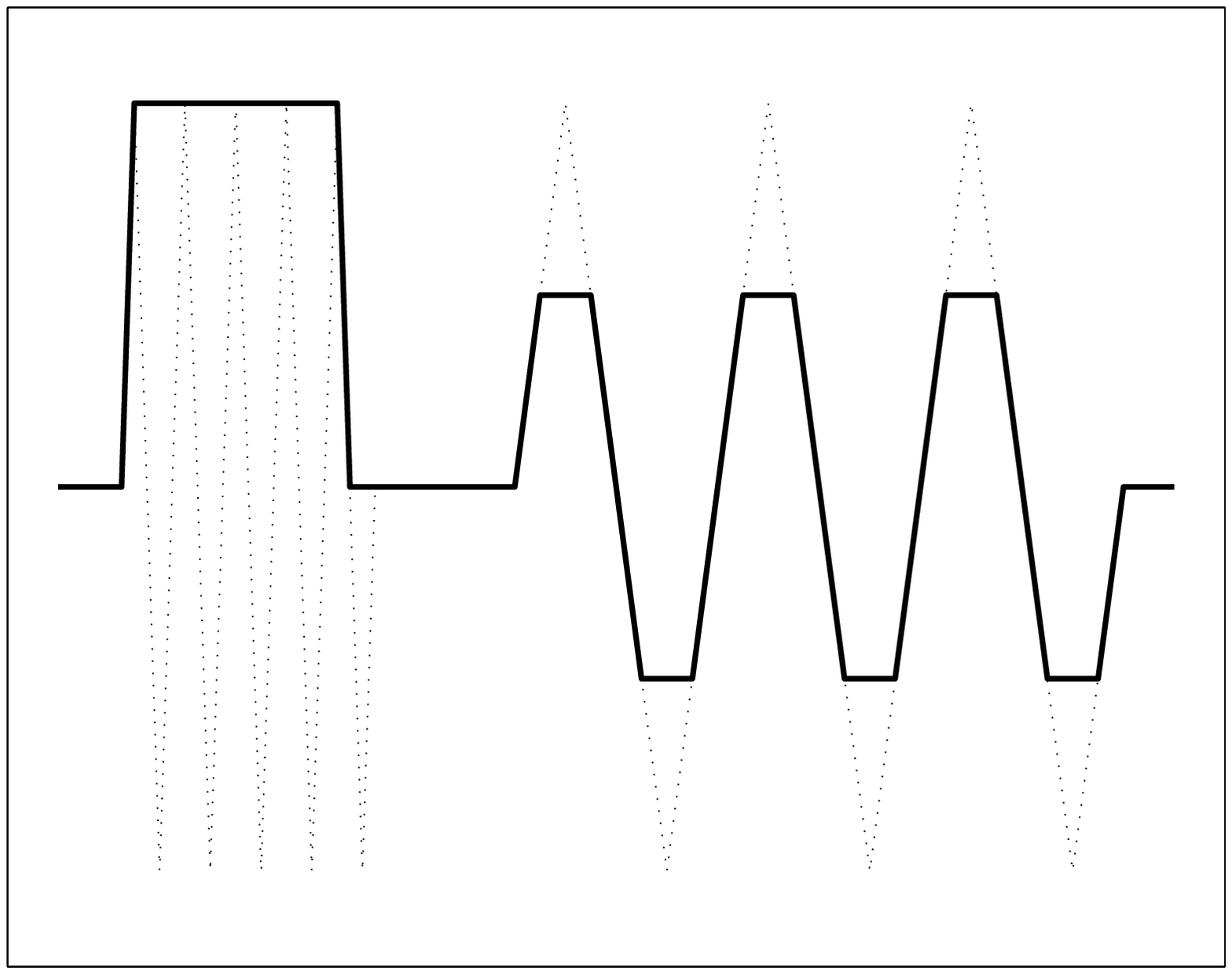}\hspace{5mm}
\includegraphics[scale=0.3]{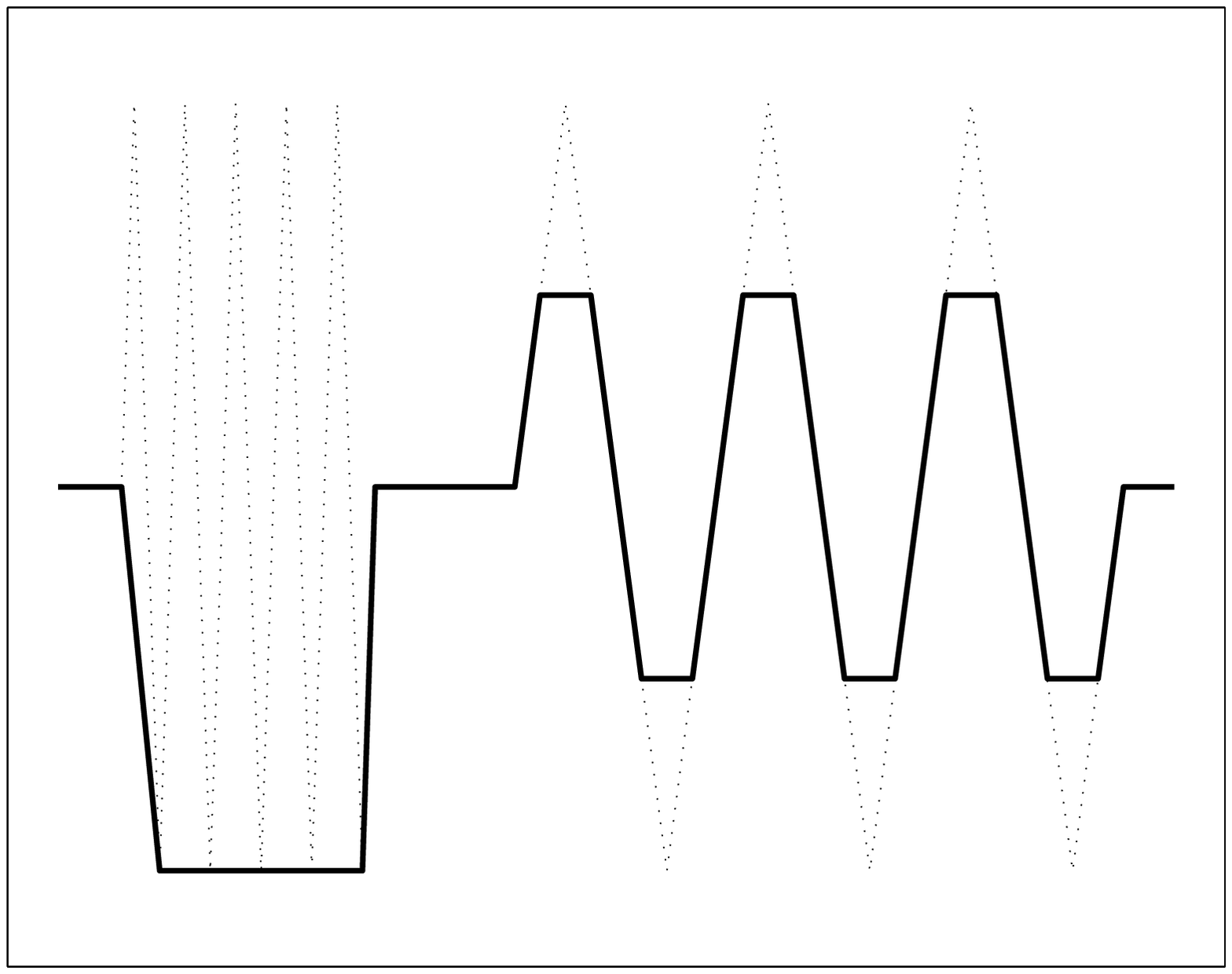}\\\nopagebreak
The functions $(L_\delta\circ U_\delta)(f)$ and $(U_\delta\circ
L_\delta)(f)$
\end{center}
\end{example}

The operator $L_\delta$ smoothes  the function $f$ from above by
removing picks while the operator $U_\delta$ smoothes the function
$f$ from below by removing pits. Note that $L_\delta\circ
U_\delta$ and $U_\delta\circ L_\delta$ resolve ambiguities in a
different way; $L_\delta\circ U_\delta$ treats oscillations of
length less then $\delta$ as picks and removes them while
$U_\delta\circ L_\delta$ considers such oscillations as pits which
are accordingly removed. The inequality $(U_\delta\circ
L_\delta)(f)\leq(L_\delta\circ U_\delta)(f)$ which is observed
here will be proved in the next section for any
$f\in\mathcal{A}(\Omega)$.

\section{The LULU semi-group}

In this section we consider the set of the operators $L_\delta$
and $U_\delta$ and their compositions. For operators on
$\mathcal{A}(\Omega)$ we consider the point-wise defined partial
order. Namely, for operators $P$, $Q$ on $\mathcal{A}(\Omega)$ we
have
\[
P\leq Q\ \Longleftrightarrow \ P(f)\leq Q(f),\
f\in\mathcal{A}(\Omega).
\]
Then the inequalities in Theorem  \ref{tLdfleqf} can be
represented in the form
\begin{equation}\label{LleqidleqU}
L_\delta\leq id\leq U_\delta,
\end{equation}
where $id$ denotes the identity operator on $\mathcal{A}(\Omega)$.

\begin{theorem}\label{tULleqLU} For any $\delta>0$ we have
$U_\delta\circ L_\delta\leq L_\delta\circ U_\delta$.
\end{theorem}
\begin{proof} Let $f\in\mathcal{A}(\Omega)$ and let $x\in\Omega$. Denote
$p=(L_\delta\circ
U_\delta)(f)(x)=S_\frac{\delta}{2}(I_\delta(S_\frac{\delta}{2}(f)))(x)$.
Let $\varepsilon$ be an arbitrary positive. For every $y\in
B_\frac{\delta}{2}(x)$ we have
\begin{equation}\label{ineqtULleqLU}
I_\delta(S_\frac{\delta}{2}(f))(y)\leq p<p+\varepsilon.
\end{equation}
\underline{Case 1.} There exists $z\in B_\frac{\delta}{2}(x)$ such
that $S_\frac{\delta}{2}(f)(z)<p+\varepsilon$. Then
$f(t)<p+\varepsilon \ \mbox{ for }\ t\in B_\frac{\delta}{2}(z)$,
which implies that $I_\frac{\delta}{2}(f)(t)< p+\varepsilon \
\mbox{ for }\ t\in B_\delta(z)$. Hence
$S_\delta(I_\frac{\delta}{2}(f))(z)\leq p+\varepsilon$. Then
$(U_\delta\circ
L_\delta)(f)(t)=I_\frac{\delta}{2}(S_\delta(I_\frac{\delta}{2}(f)))(t)\leq
p+\varepsilon$ for $t\in B_\frac{\delta}{2}(z)$. Since $x\in
B_\frac{\delta}{2}(z)$, see the case assumption, from the above
inequality we have $(U_\delta\circ L_\delta)(f)(x)\leq
p+\varepsilon$.\\
\underline{Case 2.} For every $z\in B_\frac{\delta}{2}(x)$ we have
$S_\frac{\delta}{2}(f)(z)\geq p+\varepsilon$. Denote
\[D=\left\{z\in\Omega:S_\frac{\delta}{2}(f)(z)< p+\varepsilon
\right\}.\] We will show that for every $z\in B_\delta(x)$ we have
\begin{equation}\label{intersecttULleqLU}
B_\delta(z)\cap D\neq\emptyset
\end{equation}
Due to the inequality (\ref{ineqtULleqLU}) we have that
(\ref{intersecttULleqLU}) holds for every $z\in
B_\frac{\delta}{2}(x)$. Let $z\in B_\delta(x)$ and let
$z>x+\frac{\delta}{2}$. This implies that
$x+\frac{\delta}{2}\in\Omega$. Using the inequality
(\ref{ineqtULleqLU}) for $y=x+\frac{\delta}{2}$ as well as the
case assumption we obtain that the set
$\left(x+\frac{\delta}{2},x+\frac{3\delta}{2}\right]\cap D$ is not
empty. Then $B_\delta(z)\cap D\supset
\left(x+\frac{\delta}{2},x+\frac{3\delta}{2}\right]\cap D\neq
\emptyset$. For $z<x-\frac{\delta}{2}$ condition
(\ref{intersecttULleqLU}) is proved in a similar way. Hence
(\ref{intersecttULleqLU}) holds for all $z\in B_\delta(x)$. Let
$z\in B_\delta(x)$ and $v\in B_\delta(y)\cap D$. Since $v\in D$ we
have $f(t)< p+\varepsilon$, for $t\in B_\frac{\delta}{2}(v)$.
Using that $B_\frac{\delta}{2}(z)\cap B_\frac{\delta}{2}(v)\neq
\emptyset$ we obtain that $I_\frac{\delta}{2}(f)(z)<
p+\varepsilon$, $z\in B_\delta(x)$. Therefore
$S_\delta(I_\frac{\delta}{2}(f))(x)\leq p+\varepsilon$. Then
\[
(U_\delta\circ
L_\delta)(f)(x)=I_\frac{\delta}{2}(S_\delta(I_\frac{\delta}{2}(f)))(x)\leq
S_\delta(I_\frac{\delta}{2}(f))(x)\leq p+\varepsilon.
\]
Combining the results of Case 1 and Case 2 we have $(U_\delta\circ
L_\delta)(f)(x)\leq p+\varepsilon$. Since $\varepsilon$ is
arbitrary this implies that $(U_\delta\circ L_\delta)(f)(x)\leq
p=(L_\delta\circ U_\delta)(f)(x)$.
\end{proof}

\begin{theorem}\label{tLUULidemp}
For a given $\delta>0$ the operators $L_\delta\circ U_\delta$ and
$U_\delta\circ L_\delta$ are both idempotent.
\end{theorem}
The proof is an immediate application of Lemma \ref{tIdSdId}.

\begin{theorem}\label{tULU=LU}
We have $U_\delta\circ L_\delta\circ U_\delta=L_\delta\circ
U_\delta$, $L_\delta\circ U_\delta\circ L_\delta=U_\delta\circ
L_\delta$, $\delta>0$.
\end{theorem}
\begin{proof}
Using the inequalities (\ref{LleqidleqU}) and the monotonicity of
the operators $L_\delta$, $U_\delta$, see (\ref{LdUdfmon}), we
obtain $U_\delta\circ L_\delta\circ U_\delta \geq id\circ
L_\delta\circ U_\delta\ =\ L_\delta\circ U_\delta$. For the proof
of the inverse inequality we use Theorem \ref{tULleqLU} and the
idempotence of $U_\delta$ as follows:
\[
U_\delta\circ L_\delta\circ U_\delta=(U_\delta\circ L_\delta)\circ
U_\delta\leq (L_\delta\circ U_\delta)\circ U_\delta=L_\delta\circ
(U_\delta\circ U_\delta)=L_\delta\circ U_\delta
\]
Therefore $U_\delta\circ L_\delta\circ U_\delta=L_\delta\circ
U_\delta$. The second equality is proved in a similar way.
\end{proof}

It follows from Theorems \ref{tLUULidemp} and \ref{tULU=LU} that
for a fixed $\delta>0$ every composition involving finite number
of the operators $L_\delta$ and $U_\delta$ is an element of the
set $\{L_\delta,U_\delta,U_\delta\circ L_\delta, L_\delta\circ
U_\delta\}$. Hence the operators $L_\delta$ and $U_\delta$ form a
semi-group with a composition table as follows:
\begin{center}
\begin{tabular}{c||c|c|c|c|}
&$ L_\delta$&$ U_\delta$&$U_\delta\circ L_\delta$&$L_\delta\circ
U_\delta$\\\hline\hline $L_\delta$&$ L_\delta$&$
 L_\delta\circ U_\delta$&$U_\delta\circ L_\delta$
&$ L_\delta\circ U_\delta$\\\hline $ U_\delta$&$ U_\delta\circ
L_\delta$&$ U_\delta$&$ U_\delta\circ L_\delta$ &$L_\delta\circ
U_\delta$\\\hline $ U_\delta\circ L_\delta$&$ U_\delta\circ
L_\delta$&$L_\delta\circ U_\delta$&$U_\delta\circ
L_\delta$&$L_\delta\circ U_\delta$\\\hline $ L_\delta\circ
U_\delta$&$U_\delta\circ L_\delta$&$ L_\delta\circ U_\delta$
&$U_\delta\circ L_\delta$ &$ L_\delta\circ U_\delta$\\\hline
\end{tabular}
\end{center}
Furthermore, an easy application of Theorem \ref{tULleqLU} shows
that this semi-group is completely ordered. Namely, we have
$L_\delta\leq U_\delta\circ L_\delta\leq L_\delta\circ
U_\delta\leq U_\delta$.

The smoothing of functions in $\mathcal{A}(\Omega)$ by the
compositions $ L_\delta\circ U_\delta$ and $ U_\delta\circ
L_\delta$ can be described through the concept of local
$\delta$-monotonicity.
\begin{definition}\label{deflocmon}
Let $\delta>0$. A function $f\in\mathcal{A}(\Omega)$ is called
locally $\delta$-monotone if $f$ is monotone (increasing or
decreasing) on any interval $[x,y]\subseteq \Omega$ of length not
exceeding $\delta$.
\end{definition}
\begin{theorem}\label{tlocmon}
For any given $\delta>0$ and $f\in\mathcal{A}(\Omega)$ the
functions $(L_\delta\circ U_\delta)(f)$ and $(U_\delta\circ
L_\delta)(f)$ are both locally $\delta$-monotone.
\end{theorem}
The proof uses to following technical lemma.

\begin{lemma}\label{ttechlemma02}
Let $\delta>0$ and $g\in\mathcal{A}(\Omega)$. If there exists
functions $\phi,\psi\in\mathcal{A}(\Omega)$ such that
\begin{eqnarray}
g&=&I_{\frac{\delta}{2}}(\phi),\label{ttl021}\\
g&=&S_{\frac{\delta}{2}}(\psi),\label{ttl022}
\end{eqnarray}
then $g$ is locally $\delta$-monotone.
\end{lemma}
\begin{proof}
Assume the opposite. More precisely, we assume that there exists
an interval $[a,b]\subseteq\Omega$, $b-a\leq\delta$, such that $g$
is not monotone on $[a,b]$. This means that there exists a triple
$\{x,z,y\}\subset[a,b]$, $x<z<y$, such that either
\begin{equation}\label{ttl023}
g(z)<\min\{g(x),g(y)\}\ \mbox{ or }\ g(z)>\max\{g(x),g(y)\}.
\end{equation}
Using the representation (\ref{ttl021}) and the inclusion
$B_\frac{\delta}{2}(z)\subseteq B_\frac{\delta}{2}(x)\cup
B_\frac{\delta}{2}(y)$ we obtain
\begin{eqnarray*}
g(z)&=&\sup\{\phi(t):t\in B_\frac{\delta}{2}(z)\}\ \leq\
\sup\{\phi(t):t\in B_\frac{\delta}{2}(x)\cup
B_\frac{\delta}{2}(y)\}\\
&=&\max\{\sup\{\phi(t):t\in
B_\frac{\delta}{2}(x)\},\sup\{\phi(t):t\in
B_\frac{\delta}{2}(z)\}\}\\
&=&\max\{g(x),g(y)\}.
\end{eqnarray*}
Similarly, from (\ref{ttl022}) it follows that
$g(z)\geq\min\{g(x),g(y)\}$. Thus we have
\begin{equation}\label{ttl024}
\min\{g(x),g(y)\}\leq g(z)\leq\max\{g(x),g(y)\}.
\end{equation}
The contradiction between (\ref{ttl023}) and (\ref{ttl024})
completes the proof.
\end{proof}

{\bf Proof of Theorem \ref{tlocmon}.} Let $\delta>0$ and
$f\in\mathcal{A}(\Omega)$. Denote $g=(L_\delta\circ U_\delta)(f)$.
From the definition of the operators $L_\delta$, $U_\delta$, see
(\ref{defLULU}), and Theorem \ref{tULU=LU} we obtain the following
representations of $g$:
\begin{eqnarray*}
g&=&I_\frac{\delta}{2}(S_\frac{\delta}{2}(U_\delta(f)))\\
g&=&(U_\delta\circ L_\delta\circ
U_\delta)(f)=S_\frac{\delta}{2}(I_\frac{\delta}{2}(L_\delta(U_\delta(f))))
\end{eqnarray*}
Then the local $\delta$-monotonicity of $g$ follows from Lemma
\ref{ttechlemma02}.

\section{Trend preservation}

\begin{definition}\label{defltp}
An operator $A$ is called \textbf{local trend preserving} if for
every $f\in\mathcal{A}(\Omega)$ and interval
$[x_1,x_2]\subset\Omega$ the function $A(f)$ is monotone
increasing on $[x_1,x_2]$ whenever $f$ is monotone increasing on
$[x_1,x_2]$ and $A(f)$ is monotone decreasing on $[x_1,x_2]$
whenever $f$ is monotone decreasing on $[x_1,x_2]$.
\end{definition}

\begin{definition}\label{defftp}
An operator $A$ is called \textbf{fully trend preserving} if $A$
and $id-A$ are both local trend preserving.
\end{definition}

If $A$ is local trend preserving operator then the local trend
preserving property of $id-A$ can be equivalently formulated as:
if $f$ is monotone (increasing or decreasing) on an interval
$[x_1,x_2]\subset\Omega$ then
\begin{equation}\label{defftpeq}
|A(f)(x_1)-A(f)(x_2)|\leq|f(x_1)-f(x_2)|.
\end{equation}
\begin{remark} Definition \ref{defltp} and Definition \ref{defftp} generalize the
concepts of \emph{neighbor trend preserving} and \emph{fully trend
preserving} for operators on sequences. In the context of
sequences the property (\ref{defftpeq}) is called difference
reducing, \cite{RohwerQM2002,RohwerQM2004,Rohwerbook}.
\end{remark}

\begin{theorem}\label{tcomptp}
If the operators $A$ and $B$ are fully (local) trend preserving
then so is their composition $A\circ B$.
\end{theorem}
The proof is similar to the proof of the respective statement for
sequences, see \cite[Theorem 6.10]{Rohwerbook} and will be
omitted.

We will prove that the operators $L_\delta$, $U_\delta$ and their
compositions, similar to their discrete counterparts, are all
fully trend preserving. To this end, the following technical lemma
is useful.

\begin{lemma}\label{tmon} Let function $f\in\mathcal{A}(\Omega)$ be
given and let $\delta>0$ be arbitrary.
\begin{itemize} \item[a)] If $f$ is monotone increasing on the interval $[x_1,x_2]\subseteq\Omega$ then the
function $I_\delta(f)$ is monotone increasing on
$[x_1-\delta,x_2-\delta]\cap\Omega$ and $S_\delta(f)$ is monotone
increasing on $[x_1+\delta,x_2+\delta]\cap\Omega$. \item[b)] If
$f$ is monotone decreasing on the interval
$[x_1,x_2]\subseteq\Omega$ then the function $I_\delta(f)$ is
monotone decreasing on $[x_1+\delta,x_2+\delta]\cap\Omega$ and
$S_\delta(f)$ is monotone increasing on
$[x_1-\delta,x_2-\delta]\cap\Omega$.
\end{itemize}
\end{lemma}
\begin{proof} We will prove only a) since b) is proved in a
similar way. Let $y_1,y_2\in[x_1-\delta,x_2-\delta]\cap\Omega$ and
$y_1<y_2$. We have
\begin{equation}\label{tmoneq1}
I_\delta(f)(y_1)=\inf\{f(x):x\in
[y_1-\delta,y_1+\delta]\cap\Omega\}
\end{equation}
Since $f$ is increasing on $[x_1,x_2]$ and
$[y_1+\delta,y_2+\delta]\subset[x_1,x_2]$ we have
$f(y_1+\delta)\leq f(x)$ for
$x\in[y_1+\delta,y_2+\delta]\cap\Omega$. Therefore enlarging the
interval $[y_1-\delta,y_1+\delta]$ to the interval
$[y_1-\delta,y_2+\delta]=[y_1-\delta,y_1+\delta]\cup[y_1+\delta,y_2+\delta]$
is not going to change the value of the infimum in (\ref{tmoneq1})
above. Using that the infimum of a smaller set is larger we
further have
\begin{eqnarray*}
I_\delta(f)(y_1)&=&\inf\{f(x):x\in [y_1-\delta,y_2+\delta]\cap\Omega\}\\
&\leq& \inf\{f(x):x\in
[y_2-\delta,y_2+\delta]\cap\Omega\}\\
&=&I_\delta(f)(y_2)
\end{eqnarray*}
This shows that $I_\delta(f)$ is monotone increasing on
$[x_1-\delta,x_2-\delta]\cap\Omega$. We prove that $S_\delta(f)$
is monotone increasing on $[x_1+\delta,x_2+\delta]\cap\Omega$
using a similar approach. Let
$y_1,y_2\in[x_1+\delta,x_2+\delta]\cap\Omega$ and $y_1<y_2$. By
the monotonicity of $f$ on the interval
$[y_1-\delta,y_2-\delta]\subset[x_1,x_2]$ we have
\begin{eqnarray*}
S_\delta(f)(y_2)&=&\sup\{f(x):x\in[y_2-\delta,y_2+\delta]\cap\Omega\}\\
&=&\sup\{f(x):x\in[y_1-\delta,y_2+\delta]\cap\Omega\}\\
&\geq&\sup\{f(x):x\in[y_1-\delta,y_1+\delta]\cap\Omega\}\\
&=&S_\delta(f)(y_1).
\end{eqnarray*}
\end{proof}

\begin{theorem}\label{tftp}
For an arbitrary $\delta>0$ the operators $L_\delta$, $U_\delta$
and their compositions are all fully trend preserving.
\end{theorem}
\begin{proof}
We will prove only that $L_\delta$ is fully trend preserving since
the proof of the statement for $U_\delta$ is done in a similar
way. Then, the fully trend preserving property of the compositions
follows from Theorem \ref{tcomptp}. Therefore it is sufficient to
show that if a function $f\in\mathcal{A}(\Omega)$ is monotone
increasing or monotone decreasing on an interval $[x_1,x_2]$ then
so are the functions $L_\delta(f)$ and $(id-L_\delta)(f)$. Due to
the analogy we will only discuss the situation when $f$ is
increasing.

Let $f$ be monotone increasing on $[x_1,x_2]$.

\noindent\underline{A. Proof that $L_\delta(f)$ is monotone
increasing on $[x_1,x_2]$.}

Applying Lemma \ref{tmon} a) to the operator $I_\frac{\delta}{2}$
we obtain that $I_\frac{\delta}{2}(f)$ is monotone increasing on
the interval
$[x_1-\frac{\delta}{2},x_2-\frac{\delta}{2}]\cap\Omega$.

\noindent\underline{Case 1.}
$[x_1-\frac{\delta}{2},x_2-\frac{\delta}{2}]\subset\Omega$

Using again Lemma \ref{tmon} a) for the operator
$S_\frac{\delta}{2}$ applied to $I_\frac{\delta}{2}(f)$ on the
interval $[x_1-\frac{\delta}{2},x_2-\frac{\delta}{2}]$ we obtain
that $L_\delta(f)=S_\frac{\delta}{2}(I_\frac{\delta}{2}(f))$ is
monotone increasing on $[x_1,x_2]$.

\noindent\underline{Case 2.}
$[x_1-\frac{\delta}{2},x_2-\frac{\delta}{2}]\cap\Omega=\emptyset$

Let $a$ be the left endpoint of the interval $\Omega$. For clarity
of the exposition we assume that $a\in\Omega$ but the argument
also holds if this is not true. It is easy to see that for any
$g\in\mathcal{A}(\Omega)$ the function $S_\frac{\delta}{2}(g)$ is
monotone increasing on the interval
$\left[a,a+\frac{\delta}{2}\right]$. Indeed, for
$x\in\left[a,a+\frac{\delta}{2}\right]$ we have
\[
S_\frac{\delta}{2}(g)(x)=\sup\left\{g(y):y\in\left[a,x+\frac{\delta}{2}\right]\right\}.
\]
where an increase in $x$ enlarges the interval
$\left[a,x+\frac{\delta}{2}\right]$ resulting in a higher value of
the supremum. The case assumption implies that
$[x_1,x_2]\subset\left[a,a+\frac{\delta}{2}\right]$. Since
$L_\delta(f)=S_\frac{\delta}{2}(I_\frac{\delta}{2}(f))$ is
increasing on $\left[a,a+\frac{\delta}{2}\right]$, it is also
increasing on the subinterval $[x_1,x_2]$.

\noindent\underline{Case 3.} If neither of the assumptions in Case
1 and Case 2 hold one obtains the monotonicity of $L_\delta(f)$ on
$[x_1,x_2]$ by applying Case 1 and Case 2 to suitable subintervals
of $[x_1,x_2]$.

\noindent\underline{B. Proof that $(id-L_\delta)(f)$ is monotone
increasing on $[x_1,x_2]$.}

Let $y_1,y_2\in[x_1,x_2]$, $y_1<y_2$. It follows from Part A of
the proof that
\begin{equation}\label{tftpeq0}
L_\delta(f)(y_1)\leq L_\delta(f)(y_2).
\end{equation}

\noindent\underline{Case 1.} $L(f)(y_1)=f(y_1)$. Then using that
$L_\delta(f)(y_2)\leq f(y_2)$ we obtain
\[(id-L_\delta(f))(y_1)=f(y_1)-L_\delta(f)(y_1)=0\leq f(y_2)-L_\delta(f)(y_2)=(id-L_\delta(f))(y_2)\]

\noindent\underline{Case 2.} $L(f)(y_1)<f(y_1)$. Then we have
\begin{equation}\label{tftpeq1}
I_\frac{\delta}{2}(f)(x)\leq L_\delta(f)(y_1)<f(y_1)\ \mbox{ for
all }\ x\in
\left[y_1-\frac{\delta}{2},y_1+\frac{\delta}{2}\right]\bigcap\Omega.
\end{equation}
In particular,
\begin{equation}\label{tftpeq11}
I_\frac{\delta}{2}(f)\left(y_1+\frac{\delta}{2}\right)=\inf\{f(x):
x\in[y_1,y_1+\delta]\cap\Omega\}\leq L_\delta(f)(y_1)<f(y_1).
\end{equation}
Considering the monotonicity of $f$ on the interval $[x_1,x_2]$
the above inequality implies that $y_1+\delta>y_2$. It further
follows from (\ref{tftpeq11}) that for every $\varepsilon>0$ there
exists $y_\varepsilon\in [y_2,y_1+\delta]\cap\Omega$ such that
\begin{equation}\label{tftpeq2}
f(y_\varepsilon)\leq L_\delta(f)(y_1)+\varepsilon.
\end{equation}
Hence we have
\begin{eqnarray*}
&&I_\frac{\delta}{2}(f)(x)\leq L_\delta(f)(y_1)\ ,\ \
x\in\left[y_2-\frac{\delta}{2},y_1+\frac{\delta}{2}\right]\bigcap\Omega\
\ \ \
\mbox{(see (\ref{tftpeq1}))},\\
&&I_\frac{\delta}{2}(f)(x)\leq f(y_\varepsilon)\leq
L_\delta(f)(y_1)+\varepsilon\ ,\ \ x\in
\left[y_1+\frac{\delta}{2},y_2+\frac{\delta}{2}\right]\bigcap\Omega\
\ \ \ \mbox{(see (\ref{tftpeq2}))}.
\end{eqnarray*}
Therefore
\begin{equation}\label{tftpeq3}
L_\delta(f)(y_2)=\sup\left\{I_\frac{\delta}{2}(f)(x):
x\in\left[y_2-\frac{\delta}{2},y_2+\frac{\delta}{2}\right]\bigcap\Omega\right\}\leq
L_\delta(f)(y_1)+\varepsilon.
\end{equation}
Since $\varepsilon$ in  the inequality (\ref{tftpeq3}) is
arbitrary, using also (\ref{tftpeq0}) we obtain
$L_\delta(f)(y_2)=L_\delta(f)(y_1)$. Then by the monotonicity of
$f$ on $[x_1,x_2]$ we have
\[
(id-L_\delta(f))(y_1)=f(y_1)-L_\delta(f)(y_1)\leq
f(y_2)-L_\delta(f)(y_2)=(id-L_\delta(f))(y_2).
\]

\end{proof}

\section{Total variation preservation}

The operators $L_\delta$, $U_\delta$ and their compositions are
smoothers. Therefore, one can expect that they reduce the Total
Variation of the functions. This is indeed true, but in fact these
operators satisfy a much stronger property. Namely, total
variation preservation. Denote by $BV(\Omega)$ the set of all real
functions with bounded variation defined on $\Omega$ and denote by
$TV(f)$ the total variation of $f\in BV(\Omega)$. Consider an
operator $A:BV(\Omega)\to BV(\Omega)$. Since the total variation
is a semi-norm on $BV(\Omega)$ we have
\begin{equation}\label{seminorm}
TV(f)\leq TV(A(f))+TV((id-A)(f))\ , \ \ f\in BV(\Omega).
\end{equation}

\begin{definition}\label{deftvp}
The operator $A$ is called \emph{total variation preserving} if
\begin{equation}\label{tvpreserving}
TV(f)=TV(A(f))+TV((id-A)(f))\ , \ \ f\in BV(\Omega).
\end{equation}
\end{definition}

The above definition implies that for a total variation preserving
operator the decomposition $f=A(f)+(id-A)(f)$ does not create
additional total variation.

\begin{theorem}\label{tvpcomp}
If the operators $A:BV(\Omega)\to BV(\Omega)$ and $B:BV(\Omega)\to
BV(\Omega)$ are both total variation preserving then so is their
composition $A\circ B$.
\end{theorem}
\begin{proof}
Using the total variation preserving property of $A$ and $B$ and
(\ref{seminorm}) we have
\begin{eqnarray*}
TV(f)&=&TV(B(f))+TV((id-B)(f))\\
&=&TV(A(B(f)))+TV((id-A)(B(f)))+TV((id-B)(f))\\
&\geq&TV((A\circ B)(f))+TV(((id-A)\circ B+id-B)(f))\\
&=&TV((A\circ B)(f))+TV((id-A\circ B)(f))
\end{eqnarray*}
From (\ref{seminorm}) we also obtain $TV(f)\leq TV((A\circ
B)(f))+TV((id-A\circ B)(f))$. Therefore $TV(f)=TV(A\circ
B(f))+TV((id-A\circ B)(f))$.
\end{proof}

It is easy to see that $BV(\Omega)\subseteq \mathcal{A}(\Omega)$.
Hence the operators $L_\delta$, $U_\delta$ are defined on
$BV(\Omega)$. We will show that $L_\delta$, $U_\delta$ and their
compositions are total variation preserving. The proof uses the
following technical lemmas:

\begin{lemma}\label{ttechlemma1}
Let $a,b\in\Omega$, $a\leq b$.
\begin{itemize}
\item[(a)] If there exists $\varepsilon>0$ such that
$f(x)-L_\delta(f)(x)\geq\varepsilon$, $x\in [a,b]$, then
$b-a<\delta$ and $L_\delta(f)(x)$ is a constant on $[a,b]$.
\item[(b)] If there exists $\varepsilon>0$ such that
$U_\delta(f)(x)-f(x)\geq\varepsilon$, $x\in [a,b]$, then
$b-a<\delta$ and $U_\delta(f)(x)$ is a constant on $[a,b]$.
\end{itemize}
\end{lemma}
\begin{proof} We will prove (a). Assume that $b-a\geq \delta$. Then
\[B_{\frac{\delta}{2}}\left(\frac{a+b}{2}\right)=\left[\frac{a+b-\delta}{2},\frac{a+b+\delta}{2}\right]\subseteq[a,b]
\]
and using Lemma \ref{tIdSdId} we obtain a contradiction as
follows:
\begin{eqnarray*}
I_{\frac{\delta}{2}}(f)\left(\frac{a+b}{2}\right)&=&I_{\frac{\delta}{2}}(L_\delta(f))\left(\frac{a+b}{2}\right)
=\inf_{y\in
\left[\frac{a+b-\delta}{2},\frac{a+b+\delta}{2}\right]}
L_\delta(f)(y)\\
&\leq&\inf_{y\in\left[\frac{a+b-\delta}{2},\frac{a+b+\delta}{2}\right]}f(y)-\varepsilon\
=\ I_{\frac{\delta}{2}}(f)\left(\frac{a+b}{2}\right)-\varepsilon
\end{eqnarray*}
Therefore $b-a<\delta$. Let
$p=\sup\limits_{y\in\left[b-\frac{\delta}{2},a+\frac{\delta}{2}\right]}I_\frac{\delta}{2}(f)(y)$.
Since $\left[b-\frac{\delta}{2},a+\frac{\delta}{2}\right]\subseteq
B_\frac{\delta}{2}(x)$, $x\in[a,b]$, we have
\[
p\leq\sup_{y\in
B_\frac{\delta}{2}(x)}I_\frac{\delta}{2}(f)(y)=L_\delta(f)(x)\leq
f(x)-\varepsilon\ ,\ \ x\in[a,b].
\]
Therefore
\begin{equation}\label{ttl1}
p\leq\inf_{z\in[a,b]}L_\delta(f)(z)\leq\inf_{z\in[a,b]}f(z)-\varepsilon\
,\ \ x\in[a,b].
\end{equation}
We will show next that
\begin{equation}\label{ttl2}
I_\frac{\delta}{2}(f)(y)\leq p \ \mbox{ for all } \ y\in
\left[a-\frac{\delta}{2},b+\frac{\delta}{2}\right]\cap\Omega.
\end{equation}
If $y\in\left[b-\frac{\delta}{2},a+\frac{\delta}{2}\right]$ the
inequality (\ref{ttl2}) follows directly from the definition of
$p$. Let $y>a+\frac{\delta}{2}$. Then
$[b,a+\delta]\cap\Omega\subset B_\frac{\delta}{2}(y)$ which
implies
\begin{equation}\label{ttl3}
I_\frac{\delta}{2}(f)(y)\leq\inf_{z\in[b,a+\delta]\cap\Omega}f(z).
\end{equation}
Furthermore, using (\ref{ttl1}), we have
\[
p\geq
I_\frac{\delta}{2}\left(a+\frac{\delta}{2}\right)=\min\left\{\inf_{z\in[a,b]}f(z),\inf_{z\in[b,a+\delta]\cap\Omega}f(z)\right\}
\geq\min\left\{p+\varepsilon,\inf_{z\in[b,a+\delta]\cap\Omega}f(z)\right\}.
\]
Hence
\begin{equation}\label{ttl4}
\inf_{z\in[b,a+\delta]\cap\Omega}f(z)\leq p.
\end{equation}
The inequality (\ref{ttl2}) follows from (\ref{ttl3}) and
(\ref{ttl4}). The case $y<a+\frac{\delta}{2}$ is considered in a
similar manner.

Since $B_\frac{\delta}{2}(x)\subset
\left[a-\frac{\delta}{2},b+\frac{\delta}{2}\right]$, $x\in[a,b]$,
using (\ref{ttl2}) we obtain
\begin{equation}\label{ttl5}
L_\delta{f}(x)=\sup_{y\in B_\frac{\delta}{2}(x)}
I_\frac{\delta}{2}(f)(y)\leq p\ ,\ \ x\in[a,b].
\end{equation}
The inequalities (\ref{ttl1}) and (\ref{ttl5}) imply that
$L_\delta(f)(x)=p$ for $x\in[a,b]$.
\end{proof}

\begin{lemma}\label{ttechlemma2}
Let $a,b\in\Omega$, $a\leq b$.
\begin{itemize}
\item[(a)]If $L_\delta(f)(a)\neq L_\delta(f)(b)$ then for every
$\varepsilon>0$ there exists $c\in[a,b]$ such that
\begin{eqnarray*}
(i)&&f(c)\leq \min\{f(a),f(b)\}\\
(ii)&&L_\delta(f)(c)\leq\min\{L_\delta(f)(a),L_\delta(f)(b)\}+\varepsilon\\
(iii)&&(id-L_\delta)(f)(c)\leq\min\{(id-L_\delta)(f)(a),(id-L_\delta)(f)(b)\}
\end{eqnarray*}
\item[(b)]If $U_\delta(f)(a)\neq U_\delta(f)(b)$ then there exists
$c\in[a,b]$ such that
\begin{eqnarray*}
(i)&&f(c)\geq \max\{f(a),f(b)\}\\
(ii)&&U_\delta(f)(c)\geq\max\{U_\delta(f)(a),U_\delta(f)(b)\}\\
(iii)&&(id-U_\delta)(f)(c)\geq\max\{(id-U_\delta)(f)(a),(id-U_\delta)(f)(b)\}
\end{eqnarray*}
\end{itemize}
\end{lemma}
\begin{proof}
We will prove (a) when $L_\delta(f)(a)<L_\delta(f)(b)$. The rest
is done in a similar way.

Let $\mathcal{D}=\{y\geq
a:\inf_{z\in[a,y]}(f(z)-L_\delta(f)(z))>0\}$. It follows from
Lemma \ref{ttechlemma1} that for every $y\in\mathcal{D}$ the
function $L_\delta(f)$ is a constant of $[a,y]$ and
$y-a\leq\delta$. Therefore, $b$ and $a+\delta$ are upper bounds of
$\mathcal{D}$ and we have
\[
d=\sup\mathcal{D}\leq \min\{b,a+\delta\}.
\]
Moreover, for every $\eta>0$ we have
\begin{equation}\label{ttll1}
\inf_{z\in [a,d+\eta]}(f(z)-L_\delta(f)(z))=0.
\end{equation}
\underline{Case 1.} There exists $c\in [a,d)$ such that
\[f(c)-L_\delta(f)(c)\leq
\min\{f(a)-L_\delta(f)(a),f(b)-L_\delta(f)(b)\}.\] Then (iii) is
automatically satisfied. Furthermore, (ii) holds since
$L_\delta(f)(c)=L_\delta(f)(a)<L_\delta(f)(b)$. The inequality (i)
is a consequence of (ii) and (iii). \\
\underline{Case 2.} For
every $z\in[a,d)$ we have
\begin{equation}\label{ttllcase2}
f(z)-L_\delta(f)(z)\geq\min\{f(a)-L_\delta(f)(a),f(b)-L_\delta(f)(b)\}.
\end{equation}
According to Lemma \ref{ttechlemma1} $L_\delta(f)$ is constant of
the interval $[a,d)$. Let $L_\delta(f)(x)=p$, $x\in[a,d)$. Assume
that there exists $\xi>0$ such that $\inf_{z\in[d,d+\xi]}f(z)>p$.
Lemma \ref{ttechlemma1} implies that $d+\xi<a+\delta$. Then, using
also (\ref{ttllcase2}), $\Delta=\inf_{z\in[a,d+\xi]}f(z)>p$ and we
have
\begin{eqnarray*}
p=L_\delta(f)(a)\geq
I_\frac{\delta}{2}\left(a+\frac{\delta}{2}\right)&=&\min\{\inf_{z\in[a,d+\xi])}f(z),\inf_{z\in[d+\xi,a+\delta]}f(z)\}\\
&\geq&\min\{p-\Delta,\inf_{z\in[d+\xi,a+\delta]}f(z).\}
\end{eqnarray*}
Therefore
\begin{equation}\label{ttll2}
\inf_{z\in[d+\xi,a+\delta]}f(z)\leq p.
\end{equation}
Using similar techniques as in the proof of Lemma
\ref{ttechlemma1} the inequality (\ref{ttll2}) implies that
$L_\delta(f)(x)\leq p<p+\Delta\leq f(x)$, $x\in[a,d+\xi]$, which
contradicts the definition of $d$. Hence
\begin{equation}\label{ttll3}
\inf_{z\in[d,d+\xi]}f(z)\leq p\ , \ \ \xi>0.
\end{equation}
As a consequence of the above equation we have
\begin{equation}\label{ttll4}
L_\delta(f)(d)\leq p.
\end{equation}
The function $f$, being a function of bounded variation may have
only discontinuities of first kind, that is, the left and right
limit exist at every point. Then the inequality (\ref{ttll3})
means that
\begin{equation}\label{ttll5}
f(d)\leq p \ \mbox{ or } \ f(d^+)\leq p.
\end{equation}
 The inequality (\ref{ttll1}) can be treated in a similar manner. Under the case
assumption (\ref{ttllcase2}) the inequality (\ref{ttll1}) is
equivalent to
\[
\inf_{z\in [d,d+\eta]}(f(z)-L_\delta(f)(z))=0\ ,\ \ \eta\in[
\]
which implies that
\begin{equation}\label{ttll6}
f(d)=L_\delta(f)(d) \ \mbox{ or } \ f(d^+)=L_\delta(f)(d^+).
\end{equation}
\underline{Case 2.1}\ $f(d)>p$.\ Then we also have
$L_\delta(f)(d)\leq p<f(d)$ so that (\ref{ttll5}) and
(\ref{ttll6}) imply that $f(d^+)\leq p$ and
$f(d^+)=L_\delta(f)(d^+)$. Therefore for every $\epsilon>0$ there
exists $\mu(\epsilon)>0$ such that
\[
f(z)\leq p+\epsilon\ , \ \ f(z)-L_\delta(f)(z)<\epsilon\ , \ \
z\in(d,d+\mu)\ .
\]
Let $\epsilon=\min\left\{\varepsilon,
\frac{1}{2}(f(a)-L_\delta(f)(a)),
\frac{1}{2}(f(b)-L_\delta(f)(b))\right\}$. Then any $c\in
(d,d+\mu(\epsilon))$ satisfies the conditions (i)--(iii).\\[4pt]
\underline{Case 2.2}\ $f(d)=p$.\\[3pt]
\underline{Case 2.2.1}\ $f(d)=p=L_\delta(f)(d)$.\ Then we take
$c=d$.\\[3pt]
\underline{Case 2.2.2}\ $f(d)=p>L_\delta(f)(d)$.\ Then it follows
from (\ref{ttll5}) that $f(d^+)=L_\delta(f)(d^+)$. Assume that
$L_\delta(f)(d^+)=f(d^+)>p$. Then there exists $\eta>0$ such that
$f(z)\geq L_\delta(f)(z)\geq p$, $z\in[a,d+\eta]\setminus\{d\}$
and $f(d)=p>L_\delta(f)(d)$. It is easy to see that this is
impossible. Indeed, let $L_\delta(f)(d)<m<p$. Then there exists
$y_1\in B_\frac{\delta}{2}(a)$ and $y_2\in
B_\frac{\delta}{2}(d+\eta)$ such that
$I_\frac{\delta}{2}(f)(y_1)>m$ and $I_\frac{\delta}{2}(f)(y_2)>m$.
Using also that $m$ is a lower bound of $f$ on $[a,d+\eta]$ we
obtain that $f(z)>m$, $z\in[\alpha,\beta]$ where
$\alpha=\min\{y_1,y_2\}-\frac{\delta}{2}$,
$\beta=\max\{y_1,y_2\}+\frac{\delta}{2}$. Since $d\in[\alpha,
\beta]$ and $\beta-\alpha\geq\delta$ there exists $z\in
B_\frac{\delta}{2}(d)$ such that
$B_\frac{\delta}{2}(z)\subseteq[\alpha,\beta]$. Then
$L_\delta(f)(d)\geq I_\frac{\delta}{2}(f)(z)>m$ which is a
contradiction. Therefore $L_\delta(f)(d^+)=f(d^+)\leq p$. Then the
proof proceeds as in the Case 2.1.
\end{proof}

\begin{theorem}\label{ttvp}
For an arbitrary $\delta>0$ the operators $L_\delta$, $U_\delta$
and their compositions are all total variation preserving
operators on $BV(\Omega)$.
\end{theorem}
\begin{proof}
Let $\delta>0$. We will only prove that $L_\delta$ is total
variation preserving, since the total variation preserving
property of $U_\delta$ is proved in a similar way and the
statement for the compositions follows directly from Theorem
\ref{tvpcomp}. Let $\theta>0$ and let $\{x_1,x_2, ...x_n\}$ be an
arbitrary grid of points on $\Omega$ arranged in increasing order.
We will show that there exist a finer grid $\{y_1,y_2, ...y_m\}$,
$n\leq m<2n$, such that for every $i=1,...,m-1$ we have either
\begin{eqnarray}
f(y_i)&\geq& f(y_{i+1})\nonumber\\
L_\delta(f)(y_i)+\frac{\theta}{2n}&\geq& L_\delta(f)(y_{i+1})\label{ttvp1}\\
(id-L_\delta)(f)(y_i)&\geq& (id-L_\delta)(f)(y_{i+1})\nonumber
\end{eqnarray}
or
\begin{eqnarray}
f(y_i)&\leq& f(y_{i+1})\nonumber\\
L_\delta(f)(y_i)&\leq& L_\delta(f)(y_{i+1})+\frac{\theta}{2n}\label{ttvp2}\\
(id-L_\delta)(f)(y_i)&\leq& (id-L_\delta)(f)(y_{i+1})\nonumber
\end{eqnarray}
This result is obtained from Lemma \ref{ttechlemma2} with
$\varepsilon=\frac{\theta}{n}$. If $L_\delta(f)(x_i)=
L_\delta(f)(x_{i+1})$ trivially either (\ref{ttvp1}) or
(\ref{ttvp2}) is satisfied for the points $x_i$ and $x_{i+1}$. If
$L_\delta(f)(x_i)\neq L_\delta(f)(x_{i+1})$ then according to
Lemma \ref{ttechlemma2}(a) there exists $c_i\in[x_i,x_{i+1}]$ such
that the inequalities (\ref{ttvp1}) are satisfied for the points
$x_i$ and $c_i$ and the inequalities (\ref{ttvp2}) are satisfied
for the points $c_i$ and $x_{i+1}$. Thus by including in the grid
$\{x_1,x_2, ...x_n\}$ a point $c_i$ between $x_i$ and $x_{i+1}$
for all $i$ such that $L_\delta(f)(x_i)\neq L_\delta(f)(x_{i+1})$
we obtain a finer grid $\{y_1,y_2, ...y_m\}$ satisfying either
$(\ref{ttvp1})$ or $(\ref{ttvp2})$ for every two consecutive
points. Using this property, for every $i=1,...,m-1$ we have
\begin{eqnarray*}
&&\hspace{-5mm}|f(y_i)-f(y_{i+1})|=|[L_\delta(f)(y_i)\!-\!L_\delta(f)(y_{i+1})]+
[(id\!-\!L_\delta)(f)(y_i)\!-\!(id\!-\!L_\delta)(f)(y_{i+1})]|\\
&&\ \
\geq|L_\delta(f)(y_i)\!-\!L_\delta(f)(y_{i+1})|-\frac{\theta}{n}
+|(id\!-\!L_\delta)(f)(y_i)\!-\!(id\!-\!L_\delta)(f)(y_{i+1})|
\end{eqnarray*}
Therefore
\begin{eqnarray*}
&&\hspace{-5mm}TV(f)\geq\sum_{i=1}^{m-1}|f(y_i)-f(y_{i+1})|\\
&&\geq\sum_{i=1}^{m-1}|L_\delta(f)(y_i)\!-\!L_\delta(f)(y_{i+1})|
+\sum_{i=1}^{m-1}|(id\!-\!L_\delta)(f)(y_i)\!-\!(id\!-\!L_\delta)(f)(y_{i+1})|-\theta\\
&&\geq\sum_{i=1}^{n-1}|L_\delta(f)(x_i)\!-\!L_\delta(f)(x_{i+1})|
+\sum_{i=1}^{n-1}|(id\!-\!L_\delta)(f)(x_i)\!-\!(id\!-\!L_\delta)(f)(x_{i+1})|-\theta
\end{eqnarray*}
Since the grid $\{x_1,x_2, ...x_n\}$ and the number $\theta$ are
arbitrary, the above inequality implies
\[
TV(f)\geq TV(L_\delta(f))+TV((id-L_\delta)(f)).
\]
In view of (\ref{seminorm}) this completes the proof.
\end{proof}

\section{Conclusion}
In this paper we extended the LULU operators from sequences to
real functions defined on a real interval using the lower and
upper $\delta$-envelopes of functions. The obtained structure,
although more general than the well known LULU structure of the
discrete operators, retains some of its essential properties.

Of significant importance is the link obtained between properties
of functions and sequences that are samplings of these.
Particularly, we can easily observe that if a function $f$ has a
good approximation $Af$ that is $\delta$-monotone, then a sampling
of $Af$ at a uniform sampling interval of $h$ with $h <
\frac{\delta}{n+1}$ then the sampling is $n$-monotone, and a
Discrete Pulse Transform will have no (high)-resolution components
less than $n$. Thus we may call $Af$ a ``pulse limited'' function,
in the same sense as a sequence is called ``band limited'' in the
theory of Wavelet Analysis when there are no high frequencies
present.

Since the total variation of a function is the supremum of the
total variations of all its samplings, we can derive that  the
total variation of a sequence of samplings does not exceed that of
the function. If the functions is $\delta$-monotone they are
equal, provided the sampling interval $h$ is smaller than
$\frac{\delta}{n+1}$.

This is important in image processing, where Total Variation is
used as an appropriate norm \cite{RohwerWild}. It may be
illuminating to consider that the energy reaching the ear is
appropriate as a natural norm, where the power spectrum yields
important information for economical decomposition and storage of
auditory signals.

The eye does not even see with the total illumination as norm, but
rather the measure of contrast. It is well known that we perceive
an image in the same way under different illumination intensities.
The total Variation fits naturally as the sum of the absolute
differences of intensity between neighbouring pixels. It turns out
to be the natural norm in Discrete Pulse Transforms, as they have
a naturally associated ``Parceval Identity'' which can be
considered analogous to the Parceval Identity in Wavelet and
Fourier Transforms, which is based on the energy distribution
amongst resolution levels. We thus have a Pulse Spectrum
associated with such a $LULU$-decomposition, which is useful for
thresholding decisions for economical transportation and storage
of the essentials of an image \cite{RohwerWild}.


\begin{thebibliography}{9}


\bibitem{Hamming} R.W. Hamming, Digital Filters Prentice Hall, N.J., 1956.

\bibitem{Laurie} D.P. Laurie and C.H. Rohwer, The discrete pulse transform, SIAM J. Math. Anal., {\bf 38}(3), 2007.

\bibitem{RohwerQM2002} C. H. Rohwer, Variation reduction and
$LULU$-smoothing, Quaestiones Mathematicae {\bf 25} (2002)
163--176.

\bibitem{RohwerQM2004} C. H. Rohwer, Fully trend preserving operators,
Quaestiones Mathematicae {\bf 27} (2004) 217--230.


\bibitem{Rohwerbook} C. H. Rohwer, Nonlinear Smoothers and Multiresolution
Analysis, Birkh\"{a}user, 2005.

\bibitem{RohwerWild} C.H. Rohwer and M. Wild, $LULU$ Theory, Idempotent Stack Filters, and the Mathematics of Vision of Marr, Advances in imaging and electron physics {\bf 146} (2007) 57-162.

\bibitem{Royden} H.L. Royden, Real Analysis, Macmillan, 1969.

\bibitem{Sendov} B. Sendov, Hausdorff Approximations, Kluwer,
Boston, 1990.

\bibitem{Serra} J. Serra, Image Analysis and Mathematical
Morphology, Academic Press, London, 1982.


\bibitem{Velleman} P.F. Velleman, Robust nonlinear data smoothers: definitions and recommendations, Procc. Natl. Acad. Sci. USA {\bf 74} (2), 434-436.
\end{thebibliography}
\end{document}